\pdfoutput=1
\documentclass[12pt, a4paper, openany]{amsart}

\usepackage[T2A]{fontenc}
\usepackage[utf8]{inputenc}

\usepackage[T1]{fontenc}
\usepackage{lmodern}

\usepackage[a4paper, total={7in, 9.5in}, centering]{geometry}

\usepackage{amsmath, amssymb, amsfonts, amsthm, mathtools, yhmath, mathrsfs}
\usepackage{leftindex}
\usepackage{arydshln}

\usepackage{array, longtable, hhline}

\usepackage[svgnames]{xcolor}

\usepackage[pdftex]{graphicx}
\usepackage{tikz}
\usepackage{tikz-cd}
\usepackage{quiver}
\usetikzlibrary{angles, calc, quotes, babel, arrows.meta, decorations.markings, plotmarks}
\usepgfmodule{decorations}
\usepackage{pgfplots}
\pgfplotsset{compat=newest}
\usepackage[new]{old-arrows}

\usepackage[all]{xypic}
\xyoption{rotate}
\xyoption{pdf}

\usepackage{enumerate}
\usepackage{enumitem}

\usepackage{comment}
\usepackage{pdftexcmds}
\usepackage[disable]{todonotes}
\usepackage{microtype}  

\theoremstyle{plain}
\newtheorem{theorem}{Theorem}[section]
\newtheorem*{theorem*}{Theorem}
\newtheorem{cor}[theorem]{Corollary}
\newtheorem{prop}[theorem]{Proposition}
\newtheorem{lemma}[theorem]{Lemma}

\newtheorem{convention}[theorem]{Convention}

\theoremstyle{definition}
\newtheorem{definition}[theorem]{Definition}

\newcounter{thma}
\theoremstyle{definition}

\newtheorem{remark}[theorem]{Remark}

\def\Q{{\mathbb Q}}
\def\Z{{\mathbb Z}}

\newcommand{\coker}{\mathrm{coker}}
\renewcommand{\ker}{\mathrm{ker}}

\newcommand{\sign}{\mathrm{sign}}
\newcommand{\id}{\mathrm{id}}

\newcommand{\arr}{\rightarrow}

\newcommand{\GL}{\mathit{{\rm GL}}}

\theoremstyle{plain} 
\newcommand{\thistheoremname}{}
\newtheorem*{genericthm}{\thistheoremname} 
\newenvironment{namedthm}[1]
  {\renewcommand{\thistheoremname}{#1}%
   \begin{genericthm}}
  {\end{genericthm}}

\title{Unordered resolutions and homological stability for linear groups}
\author{Ivan Vasilev, Serge Yagunov}
\date{\today}

\begin{document}

\begin{abstract}
	In this paper, we develop a modified proof strategy for homological stability of linear groups, with the general linear groups serving as a primary example. Our arguments are more direct than those in the classical works of Quillen and Suslin--Nesterenko, although they apply only with localized coefficients. The localization at \((n-1)!\) that arises in our approach appears to be closely related to several conjectures of Mirzaii 
	as well as to Suslin’s injectivity conjecture.
\end{abstract}

\maketitle

\section{Introduction}
 
Homological stability phenomena play a central role in many areas of mathematics, ranging from algebraic $K$-theory and algebraic topology to more recent applications in number theory and geometric group theory. Notable examples include the work of Miller--Patzt--Petersen--Randal-Williams on $L$-functions~\cite{MillerPatztPetersenRandalWilliams2025} and that of Szymik--Wahl on the Higman--Thompson groups~\cite{SzymikWahl2019}.

Roughly speaking, a sequence of groups
\[
G_1 \;\hookrightarrow\; G_2 \;\hookrightarrow\; \cdots
\]
is said to \emph{exhibit homological stability} if the induced maps in homology
\[
H_i(G_n; \mathbb{Z}) \;\longrightarrow\; H_i(G_{n+1}; \mathbb{Z})
\]
are isomorphisms for all sufficiently large $n$ (depending on $i$).

One of the most fundamental and historically significant examples of such a family is given by the general linear groups $\GL_n(R)$. The study of their homology is deeply intertwined with algebraic $K$-theory and has been a driving force behind the development of the subject.

Indeed, Quillen’s definition of higher $K$-groups associates to a ring $R$ the homotopy groups of a space constructed from the infinite general linear group
\[
\GL(R) = \varinjlim_n \GL_n(R).
\]
The \emph{Hurewicz homomorphism} provides a natural bridge between homotopy and homology, yielding canonical maps
\[
h_n \colon K_n(R) \;\longrightarrow\; H_n(\GL(R); \mathbb{Z}),
\]
which relate algebraic $K$-theory to the more computable structure of group homology.

The fact that general linear groups over fields exhibit homological stability was first established by Quillen in his seminal work~\cite{quillen1972}, where he proved stability for general linear groups over finite fields. This result was subsequently generalized by van der Kallen~\cite{vdkallen1980} and by Nesterenko and Suslin~\cite{nes-sus1990}.

In the case of infinite fields, the optimal stability range known to date can be stated as follows:


\begin{theorem*}
Let \( F \) be an infinite field. Then the inclusion \( \GL_{n-1}(F) \hookrightarrow \GL_{n}(F) \) induces an isomorphism
\[
H_i(\GL_{n-1}(F), \mathbb{Z}) \;\xrightarrow{\;\cong\;}\; H_i(\GL_{n}(F), \mathbb{Z})
\]
for all \( n > i \).
\end{theorem*}

Nesterenko and Suslin showed that this stability range is sharp by explicitly describing the cokernel of the first non-isomorphic stabilization map.

\begin{theorem*}
	Let $F$ be an infinite field. Then there is a natural exact sequence
	\[
	H_{n}(\GL_{n-1}(F)) \longrightarrow H_{n}(\GL_{n}(F)) \twoheadrightarrow K^{M}_{n}(F),
	\]
	where $K^{M}_{n}(F)$ denotes the $n$-th Milnor $K$-group of $F$.
\end{theorem*}

The usual strategy for proving homological stability results can be briefly described as leveraging spectral sequences arising from group actions on suitable spaces or complexes. A comprehensive exposition of this strategy, together with numerous applications, can be found in~\cite{wahl2022}.

In this paper, we work with a modified version of the general position complex, which is standard in this area. Let $A_k$ denote the free abelian group generated by $k$-tuples of vectors in $F^{n}$ such that every subset of $\min(n,k)$ vectors is linearly independent. These free abelian groups form a complex of $\GL_{n}(F)$-modules via left multiplication, equipped with the standard differential defined by the alternating sum of face maps.

Each free abelian group $A_k$ can also be endowed with a right action of the symmetric group $S_k$, defined on generators by
\[
(v_1, \ldots, v_k) \cdot \sigma
= (-1)^{\mathrm{sign}(\sigma)}(v_{\sigma^{-1}(1)}, \ldots, v_{\sigma^{-1}(k)}),
\]
which should be viewed as an analogue of orientation. It is straightforward to verify that this action commutes with the differential, allowing us to define the quotient complex
\[
\widetilde{A}_{\ast} = (A_{\ast})_{S_\ast}.
\]

The main result of this paper is that this complex completely describes the relative homology of general linear groups up to a finite torsion part. More precisely, using our approach we reprove the following.

\begin{namedthm}{Theorem}
	For any $n \ge 1$ and any infinite field $F$ denote by $A^{(n)}_k$ the complex generated by n-type general position vectors (see below). Then there exists an isomorphism
	\[
	H_\ast(\GL_n(F), \GL_{n-1}(F); \mathbb{Q})
	\cong
	H_{\ast+1}\bigl((\widetilde{A}^{(n)}_{\ast})_{\GL_n}\bigr) \otimes \mathbb{Q}.
	\]
\end{namedthm}

When formulated more explicitly, this theorem yields two important corollaries.

\begin{cor}
	Let $F$ be an infinite field. Then the inclusion
	$\GL_{n-1}(F) \hookrightarrow \GL_{n}(F)$
	induces an isomorphism
	\[
	H_i(\GL_{n-1}(F); \mathbb{Z}[\tfrac{1}{(n-1)!}])
	\;\xrightarrow{\;\cong\;}\;
	H_i(\GL_{n}(F); \mathbb{Z}[\tfrac{1}{(n-1)!}])
	\]
	for all $i \le n - 1$.
\end{cor}

\begin{cor}\label{maincor2}
	Let $F$ be an infinite field. Then there is a natural exact sequence
	\[
	H_{n}(\GL_{n-1}(F); \mathbb{Z}[\tfrac{1}{(n-1)!}])
	\longrightarrow
	H_{n}(\GL_{n}(F); \mathbb{Z}[\tfrac{1}{(n-1)!}])
	\twoheadrightarrow
	K^{M}_{n}(F) \otimes \mathbb{Z}[\tfrac{1}{(n-1)!}],
	\]
	where $K^{M}_{n}(F)$ denotes the $n$-th Milnor $K$-group of $F$.
\end{cor}

These results coincide with the theorems of Suslin and Nesterenko. Using unordered resolutions, we obtain shorter and more explicit proofs, at the cost of localizing coefficients. The use of unordered resolutions elliminates part of the homology coming from the monomial subgroups of $\GL$.
This method is inspired by ideas developed by the second author
in his study of the homology of bi-Grassmannian complexes~\cite{yagunov1997}.
Similar technique can be applied to the case of groups $\GL_\ast(R)$, assuming that $R$ is a local ring with many units.

\medskip

\noindent\textbf{Acknowledgments.}
The first author is deeply grateful to Vasily Ionin, who significantly influenced his mathematical interests and to whom he attributes his current passion for homotopy theory and related fields. He also thanks Viktor Lavrukhin, Ivan Vorobyev, Roman Mikhailov, and the entire collective of Room~209 for providing an atmosphere conducive to mathematical insights.

I.V.\ was supported by a scholarship from the Rodnye Goroda Foundation and is currently supported by the Faculty of Mathematics at the Higher School of Economics.


\todo[inline]{Add part about Suslin's Injectivity Conjecture}

\section{Affine resolutions and their unordered versions}

\begin{convention}
Throughout this paper, we (usually) work with the general linear group over a fixed infinite field \( F \).  
For brevity, we will often denote \( \GL_{n}(F) \) by \( \GL_n \) or simply \( G_n \).
\end{convention}

Let $A_k$ denote the free abelian group generated by the set of all $n \times k$ matrices over the field $F$, i.e.
\[
A_k \coloneq \mathbb{Z}\big[M_{n \times k}(F)\big].
\]

By the convention, we set \( A_0 = \mathbb{Z} \).  
This graded abelian group is equipped with a differential defined as follows:
\begin{equation*}
  \partial_k([X_1, X_2, \ldots, X_k]) 
  = [X_2, \ldots, X_k] - [X_1, X_3, \ldots, X_k] + \cdots + (-1)^{k+1}[X_1, \ldots, X_{k-1}],
\end{equation*}
where \( X_i \) are the columns of the generator.

The group \( G_n \) acts via left multiplication on these abelian groups.  
A standard verification shows that the graded abelian group \( A_\ast \), with the differential above, forms a complex of \( G_n \)-modules.
Certain subcomplexes of this complex arise naturally in the context of homological stability. The generators of the 
\( G_n \)-modules appearing in these subcomplexes typically satisfy so-called \emph{general position conditions}, as introduced and formally defined in~\cite{Kerz2005}.

We now introduce the notion of an \emph{$i$-type general position condition} for a collection of vectors. Variants of such conditions have appeared previously with important applications; see, for example,~\cite{nes-sus1990},~\cite{sus1990}, and~\cite{mirzaii2005}. In particular, these works make essential use of what we shall refer to as the $2$-type and $n$-type general position conditions.

For a local ring \( R \):
\begin{itemize}
    \item Let \( R^\infty \) denote the \emph{free (left) \( R \)-module} with basis \( e_1, e_2, \ldots \).
    \item Let \( R^n \) denote the submodule of \( R^\infty \) with basis \( e_1, \ldots, e_n \).
    \item A vector \( v \in R^n \) is identified with its corresponding column vector of height \( n \).
\end{itemize}

Vectors \( v_1, \ldots, v_m \in R^\infty \) are called \emph{jointly unimodular} (or said to form a \emph{unimodular frame}) if they constitute a basis of a free direct summand of \( R^\infty \).  
For vectors \( v_1, \ldots, v_m \in R^n \), the property of being jointly unimodular is equivalent to the \emph{left invertibility} of the matrix \( (v_1, \ldots, v_m) \in \mathrm{M}_{n \times m}(R) \).

\begin{definition}\label{def2.2}
  A set of \( k \) vectors in \( R^{n} \) is said to satisfy the \emph{$i$-type general position condition} if every subset of \( \min(i, k) \) vectors forms a unimodular frame.
\end{definition}

\begin{remark}
	Specializing to the case where the ring $R$ is a field, a collection of vectors forms a unimodular frame if and only if the vectors are linearly independent. Although our definitions apply to arbitrary local rings, in this paper we restrict attention, for simplicity, to infinite fields.
\end{remark}

This definition allows us to define a natural filtration of the complex \( A_{\ast} \).  
Let \( A^{(i)}_{k} \) denote the free abelian group generated by all \( n \times k \) matrices whose columns satisfy the \( i \)-type general position condition.  
These abelian groups form a complex of \( G_n \)-modules \( A_{\ast}^{(i)} \), since left multiplication by an element of \( \GL_n(F) \) preserves the \( i \)-type general position condition for every \( i \).  
The resulting filtration is:
\[
A_{\ast} = A_{\ast}^{(0)} \supset A_{\ast}^{(1)} \supset \cdots \supset A_{\ast}^{(n)} \supset 0.
\]

\todo[inline]{Try to do something with this filtration}

In fact, all of these complexes are acyclic. We will provide an explicit proof for $A_{\ast}^{(n)}$, as this case constitutes the primary focus of the paper; analogous arguments apply equally to all other cases. We thus formulate the following general statement.

\begin{prop}\label{prop2.3}
For any \( 0 \leq i \leq n \), the complex \( A_{\ast}^{(i)} \) is acyclic.
\end{prop}

\begin{proof} We consider the case $i=n$.
To prove acyclicity, we explicitely demonstrate that every cycle is a boundary.  
Let
\[
c = \sum_{j} n_{j} (a_1^{j}, \ldots, a_{k}^{j}) \in A^{(n)}_{k},
\]
$k<n$, be a cycle.  

Since the field $F$ is infinite, there exists a vector $\tilde{a}$ such that, for every summand, the set of vectors
\[
(\tilde{a}, a_1^{j}, \ldots, a_{k}^{j})
\]
still forms a unimodular frame. 
Define
\[
b = \sum_{j} n_{j} (\tilde{a}, a_1^{j}, \ldots, a_{k}^{j}) \in A^{(n)}_{k+1}.
\]
One can easily verify that $c = \partial_{k+1}(b)$. 
Hence, $c$ is a boundary, and the complex $A_{\ast}^{(n)}$ is acyclic.
\end{proof}
	
We are now prepared to define the \emph{unordered} versions of the complexes introduced above.  
For each \( i \), the module \( A_{k}^{(i)} \) carries a natural right \( S_{k} \)-module structure.  
It suffices to specify the action on generators: for a generator \([X_1, X_2, \ldots, X_{k}]\) and any permutation \( \sigma \in S_{k} \), define
\[
[X_1, X_2, \ldots, X_{k}] \sigma 
  = (-1)^{\mathrm{sign}(\sigma)}[X_{\sigma^{-1}(1)}, X_{\sigma^{-1}(2)}, \ldots, X_{\sigma^{-1}(k)}].
\]

We then define the quotient complex
\[
\widetilde{A}_{k}^{(i)} \coloneqq \bigl(A^{(i)}_{k}\bigr)_{S_k}
  \;\cong\;
  A^{(i)}_{k} \bigg/ 
  \left\langle m - (-1)^{\mathrm{sign}(\sigma)}m\sigma 
  \,\bigg|\,
  m \in A^{(i)}_{k},\; \sigma \in S_{k} \right\rangle.
\]

The resulting groups form a complex of $\GL_{n}$-modules, which we refer to as the \emph{unordered resolutions}, in analogy with unordered configuration spaces. These relations can be interpreted as encoding the orientation of each tuple of vectors within the complex.

\todo[inline]{Can it be generalized?}

\begin{remark}
	These are indeed complexes of $G_n$-modules because the left $G_n$-action is compatible with the right $S_k$-action.
\end{remark}

It turns out that the unordered complexes are also acyclic, essentially by adapting the argument of Proposition~\ref{prop2.3}.

\begin{prop}\label{cor2.6}
	For every $0 \le i \le n$, the complex $\widetilde{A}_{\ast}^{(i)}$ is acyclic.
\end{prop}

The resulting groups form a complex of \( \GL_{n} \)-modules, which we refer to as the \emph{unordered resolutions}, in analogy with unordered configuration spaces.  
These relations may be interpreted as encoding the orientation of each tuple of vectors within the complex.

\todo[inline]{Can it be generalized?}

\section{Quillen’s stability argument}
Quillen developed an approach for establishing homological stability by leveraging spectral sequences arising from group actions on appropriate spaces or complexes.  
To implement this method for a family of groups $\{G_n\}_{n \ge 0}$, one constructs a corresponding family of chain complexes $\{A^{n}_{\ast}\}_{n \in \mathbb{N}}$ of $G_n$-modules that satisfy suitable acyclicity and functoriality properties (see~\cite{wahl2022} for a detailed discussion).

Let us now illustrate this approach in our particular setting. Given a complex $C_{\ast}$ of $G_n$-modules, consider the standard spectral sequence (see Chapter 7.5 of~\cite{brown2012}): 
\[
\leftindex^I{E}_{p,q}^{1} = H_{q}(G_n; C_{p}) \implies H_{p+q}(G_n; C_{\ast}).
\]  
If the complex $C_{\ast}$ is acyclic, this spectral sequence converges to zero by the following well-known proposition:  

\begin{prop}[{\cite[Proposition~5.2]{brown2012}}]
  If $\tau: C \arr C'$ is a weak equivalence of $G$-chain complexes, then $\tau$ induces an isomorphism $H_{\ast}(G, C) \cong H_{\ast}(G, C')$
\end{prop}

We work with the spectral sequences arising from the unordered resolutions introduced above. For each $i$, this yields a spectral sequence
\[
\leftindex^I{E}_{p,q}^{1} = H_{q}(G_n; \widetilde{A}_{p}^{(i)}) \implies 0.
\]  
In the next section, we analyze special properties of this spectral sequence for $i = n$. Notably, the $\leftindex^I{E}$-page contains numerous vanishing terms when the resolution is localized, which enables a streamlined proof of homological stability for the general linear group.

\section{Properties of the spectral sequences and homological stability proof}
This section contains the key statements necessary for understanding the properties of the spectral sequence arising from $\widetilde{A}_{\ast}^{(n)}$. We begin with foundational propositions from~\cite{nes-sus1990}, culminating in the key Lemma~\ref{keylemma}.

Let $M_{n \times m}(F)$ denote the set of all $n \times m$ matrices over $F$, and define the \emph{affine group} as:  
\[
\mathrm{Aff}_{n,m} \coloneqq \left\{ 
\begin{pmatrix} 
I_n & M_{n \times m}(F) \\ 
0 & \GL_m(F) 
\end{pmatrix} \right\},
\]
where $I_n$ is the $n \times n$ identity matrix.

\begin{theorem}[{\cite[Theorem~1.11]{nes-sus1990}}]\label{affine_iso}
  The inclusion $\GL_n(F) \hookrightarrow \mathrm{Aff}_{m,n}$ induces an isomorphism
  \[
  H_{\ast}(\GL_n(F); \mathbb{Z}) \cong H_{\ast}(\mathrm{Aff}_{m,n}; \mathbb{Z}).
  \]
\end{theorem}

The homology of affine groups is particularly significant, as it appears naturally on the $E^1$-page of our spectral sequence.  
We present the following well-known lemma (see~\cite{nes-sus1990}) together with a proof, since the technique will recur in the subsequent discussion.  
Before stating it, however, we first require the following structural result.

\begin{prop}[{\cite[Corollary~5.4]{brown2012}}]\label{prop4.2}
Let \( N \) be a \( G \)-module whose underlying abelian group is of the form \(\bigoplus_{i \in I} M_i\). Assume that the \( G \)-action permutes the summands according to some action of \( G \) on \( I \). Let \( G_i \) be the isotropy group of \( i \) and let \( E \) be a set of representatives for 
\(I\bmod G\). Then \( M_i \) is a \( G_i \)-module and there is a \( G \)-isomorphism 
\[
N \cong \bigoplus_{i \in E} \operatorname{Ind}_{G_i}^G M_i.
\]
\end{prop}

\begin{lemma}[{\cite[Lemma~2.3]{nes-sus1990}}]\label{lemma4.2}
  There is an isomorphism $H_{*}(\GL_n; \widetilde{A}^{(n)}_{1}) \cong H_{*}(\GL_{n-1}; \Z)$.  
  The map  
  \[
  H_{*}(\GL_{n-1}; \Z) \xrightarrow{\sim} H_{*}(\GL_n; \widetilde{A}^{(n)}_{1}) \xrightarrow{\epsilon} H_{*}(\GL_n; \Z)
  \]  
  is induced by the standard embedding $\GL_{n-1} \hookrightarrow \GL_n$.\qed
\end{lemma}

\begin{proof}
First observe that $\widetilde{A}_{1}^{(i)} = A_{1}^{(i)}$ for all $i$, since $S_{1}$ is trivial. In particular, $\widetilde{A}_{1}^{(n)}$ coincides with Nesterenko--Suslin's $C_1$. The group $\GL_n$ acts transitively on $\widetilde{A}_{1}^{(n)}$, and by combining Shapiro's lemma with Proposition~\ref{prop4.2}, we obtain
\[
H_{\ast}(\GL_n, \widetilde{A}_{1}^{(n)}) \cong H_{\ast}(\mathrm{Stab}_{\GL_n}(e_1), \Z),
\]
where $e_1 = (1, 0, \ldots, 0)^\intercal$. The stabilizer $\mathrm{Stab}_{\GL_n}(e_1)$ is precisely the affine group $\mathrm{Aff}_{1,n-1}$. Applying Theorem~\ref{affine_iso} then yields
\[
H_{\ast}(\mathrm{Aff}_{1,n-1}; \Z) \cong H_{\ast}(\GL_{n-1}; \Z),
\]
establishing thus the first isomorphism. The second claim follows directly from \cite[Lemma~2.4]{nes-sus1990}.
\end{proof}

\begin{remark}
  Similarly, one establishes a stronger result~\cite{nes-sus1990}:  
  \[
  H_{*}(\GL_{n}; A_{k}^{(n)}) \cong H_{*}(\GL_{n-k}; \mathbb{Z})
  \]
  The proof uses arguments analogous to those in Lemma~\ref{lemma4.2}. We present this as a remark, since our treatment primarily focuses on the unordered variant of the complexes.

\end{remark}

Now we are ready to prove the key lemma:

\begin{lemma}\label{keylemma}
  The map $(\id, \sigma)\colon H_{*}(\GL_n; A^{(n)}_{k}) \rightarrow H_{*}(\GL_n; A^{(n)}_{k})$ is the identity for all $2 \leq k \leq n$ and any permutation $\sigma \in S_{k}$.
\end{lemma}
\begin{proof}
  It suffices to prove the statement for a transposition \(\theta\) of two arbitrary adjacent columns.  
 Let $u = (e_{n - k + 1}, \ldots, e_n)$.  
There exists an element $h \in \GL_n$ --- namely, the matrix whose left multiplication swaps two adjacent rows of $u$ --- which consequently swaps the columns corresponding to $\theta$.  
This is possible because $u$ consists of the last $k$ columns of the identity matrix.

  Consider the following commutative diagram in the category of pairs (group, module):  
  
  \[\begin{tikzcd}
    {(\GL_{n-k}, \Z)} &&& {(\GL_{n}, A^{(n)}_{k})} \\
    \\
    &&& {(\GL_{n}, A^{(n)}_{k})}
    \arrow["{(i, u)}", from=1-1, to=1-4]
    \arrow[from=1-1, to=3-4]
    \arrow["{(hgh^{-1}, hm)}", shift left=2, from=1-4, to=3-4]
    \arrow["{(id, \sigma)}"', shift right=4, from=1-4, to=3-4]
  \end{tikzcd}\]

  The top arrow is induced by the standard embedding $\GL_{n-k} \hookrightarrow \GL_n$ and the map sending $1 \mapsto u$ on modules.  
  The composition of the top arrow with the two vertical arrows coincides, which forces the diagonal arrow to be the unique map making the diagram commute.  
  
  Applying the homology functor then yields:

  \[\begin{tikzcd}
	{H_{*}(\GL_{n-k}, \Z)} &&& {H_{*}(\GL_{n}, A^{(n)}_{k})} \\
	\\
	&&& {H_{*}(\GL_{n}, A^{(n)}_{k})}
	\arrow["{(i, u)_{*}}", from=1-1, to=1-4]
	\arrow[from=1-1, to=3-4]
	\arrow["{(hgh^{-1}, hm)_{*}}", shift left=2, from=1-4, to=3-4]
	\arrow["{(id, \sigma)_{*}}"', shift right=4, from=1-4, to=3-4]
\end{tikzcd}\]
   
  By \cite[§3.8]{brown2012}, the map $(hgh^{-1}, hm)_{*}$ is the identity on homology.  
   Lemma \ref{lemma4.2} implies $(i, u)_{*}$ is an isomorphism.  
  Commutativity forces the diagonal arrow to be an isomorphism, hence $(\id, \sigma)_{*}$ acts identically on homology.  
\end{proof}

The preceding lemma plays a key role in analyzing the structure of our spectral sequence, as illustrated by the following results.

\begin{lemma}\label{lemma4.4}
    Let \( 2 \leq k \leq n \). Then \( H_{*}(\GL_n, \widetilde{A}^{(n)}_k) \otimes \Z[\tfrac{1}{k!}] = 0 \). 
\end{lemma}

\begin{proof}
    Let \(\sigma_i\) denote the transposition swapping the \(i\)-th and \((i+1)\)-th columns. We work with the localized complex where \(k!\) is inverted in \(A^{(n)}_{*}\).  

    The module \(\widetilde{A}^{(n)}_k\) is realized as the cokernel of:  
    \[
    \bigoplus_{i=1}^{k-1} A^{(n)}_{k} \xrightarrow{\phi} A^{(n)}_{k} \twoheadrightarrow \widetilde{A}^{(n)}_k \to 0,
    \]  
    where \(\phi = \sum_{i=1}^{k-1} (\id + \sigma_i)\) acts via \((m_1, \ldots, m_{k-1}) \mapsto \sum_{i=1}^{k-1}(m_i + m_i\sigma_i)\).  

    The cokernel captures precisely the antisymmetry relations
    \[
    m = (-1)^{\mathrm{sign}(\sigma)} m \sigma,
    \]
    generated by transpositions. The long exact homology sequence associated to the short exact sequence
    \[
    0 \longrightarrow \mathrm{Im}(\phi) \hookrightarrow A^{(n)}_{k} \twoheadrightarrow \widetilde{A}^{(n)}_k \longrightarrow 0
    \]
    yields:
     \[
    H_{*}(\GL_n; \mathrm{Im}(\phi)) \xrightarrow{j_{*}} H_{*}(\GL_n; A^{(n)}_{k}) \xrightarrow{p_{*}} H_{*}(\GL_n; \widetilde{A}^{(n)}_k).
    \]  
    We prove that \(j_*\) becomes an isomorphism after localization.

    \textbf{Surjectivity}.  
    The composition:  
    \[
    H_{*}\left(\GL_n; \bigoplus_{i=1}^{k-1} A^{(n)}_{k}\right) \xrightarrow{\phi_*} H_{*}(\GL_n; \mathrm{Im}(\phi)) \xrightarrow{j_*} H_{*}(\GL_n; A^{(n)}_{k})
    \]  
    satisfies \(\phi_* = \bigoplus_{i=1}^{k-1} (\id + \id) = 2 \cdot \bigoplus \id\) by Lemma~\ref{keylemma}.  Since 2 is invertible in the coefficint ring,
    the surjectivity is proven.

    \textbf{Injectivity}.  
    Define \(\psi: A^{(n)}_{k} \to \mathrm{Im}(\phi)\) by:  
    \[
    \psi \coloneqq \sum_{\sigma \in S_k} (\id + (-1)^{\mathrm{sign}(\sigma)}\sigma).
    \]  
    The composition:  
    \[
    H_{*}(\GL_n; \mathrm{Im}(\phi)) \xrightarrow{j_*} H_{*}(\GL_n; A^{(n)}_{k}) \xrightarrow{\psi_*} H_{*}(\GL_n; \mathrm{Im}(\phi))
    \]  
    acts on generators as:  
    \begin{align*}
        \psi(m + m\sigma_i) &= \sum_{\sigma \in S_k} \left[(\id - (-1)^{\mathrm{sign}(\sigma)}\sigma)(m + m\sigma_i)\right] \\
        &= k!(m + m\sigma_i) - \sum_{\sigma \in S_k} (-1)^{\mathrm{sign}(\sigma)}m\sigma - \sum_{\sigma \in S_k} (-1)^{\mathrm{sign}(\sigma)}m\sigma_i\sigma \\
        &= k!(m + m\sigma_i).
    \end{align*}  
    Since any transvection is odd, the sign-alternating terms are annihilated, leaving multiplication by $k!$, which becomes an isomorphism after localization. Hence, $j_*$ is injective. 
\end{proof}

\todo[inline]{Is it true because of the statement like "homology groups of $\GL_n$ with coefficients in a functor of a finite degree vanish rationally"?}

\begin{lemma}\label{lemma4.6}
  For any \( k \geq n + 1 \),
  \[
  H_{i}(\GL_n; \widetilde{A}^{(n)}_k) \otimes \Z[\tfrac{1}{k!}] = 
  \begin{cases}
  H_{i+1}(\widetilde{A}^{(n)}_k)_{\GL_n} \otimes \Z[\tfrac{1}{k!}], & i = 0, \\
      0, & \text{otherwise}.
  \end{cases}
  \]
\end{lemma}

\begin{proof}
First observe that \(\GL_n\) acts freely on the generators of \(A^{(n)}_{k}\) for \(k \geq n + 1\), making \(A^{(n)}_{k}\) a free \(\Z[\GL_n]\)-module by Proposition \ref{prop4.2}, hence \(H_{i \geq 1}(\GL_n; A^{(n)}_{k}) = 0\). Under the assumption that \(k!\) is inverted in both \(\widetilde{A}^{(n)}_{*}\) and \(A^{(n)}_{*}\), we establish:
\[
H_{i}(\GL_n; \widetilde{A}^{(n)}_k) = 
\begin{cases}
    H_{i+1}(\widetilde{A}^{(n)}_k)_{\GL_n}, & i = 0, \\
    0, & \text{otherwise}.
\end{cases}
\]

Consider the long exact sequence from Lemma~\ref{lemma4.4}:
\begin{equation}\label{les_1}
    \cdots \to H_{*}(\GL_n; \mathrm{Im}(\phi)) \xrightarrow{j_{*}} H_{*}(\GL_n; A^{(n)}_{k}) \xrightarrow{p_{*}} H_{*}(\GL_n; \widetilde{A}^{(n)}_k) \to \cdots
\end{equation}
To show \(H_{i \geq 2}(\GL_n; \widetilde{A}^{(n)}_k) = 0\), we first prove \(H_{i \geq 1}(\GL_n; \mathrm{Im}(\phi)) = 0\). From Lemma~\ref{lemma4.4}, the composition:
\begin{equation}\label{les_2}
    H_{*}(\GL_n; \mathrm{Im}(\phi)) \xrightarrow{j_{*}} H_{*}(\GL_n; A^{(n)}_{k}) \xrightarrow{\psi_\ast} H_{*}(\GL_n; \mathrm{Im}(\phi))
\end{equation}
is an isomorphism (multiplication by \(k!\)). Since \(H_{i \geq 1}(\GL_n; A^{(n)}_{k}) = 0\), this implies \(H_{i \geq 1}(\GL_n; \mathrm{Im}(\phi)) = 0\).

For \(H_1(\GL_n; \widetilde{A}^{(n)}_k) = 0\), we examine the residual segment of (\ref{les_1}):
\[
0 \to H_{1}(\GL_n; \widetilde{A}^{(n)}_k) \to (\mathrm{Im}(\phi))_{\GL_n} \xrightarrow{j_{*}} (A^{(n)}_{k})_{\GL_n} \to (\widetilde{A}^{(n)}_k)_{\GL_n} \to 0.
\]
Since $j_*$ is injective (by (\ref{les_2})), we conclude that $H_1(\GL_n; \widetilde{A}^{(n)}_k) = 0$, as required.
\end{proof}





The preceding lemmas suffice to establish one of the central results of this paper.  
We now show that the complex $\widetilde{A}_{\ast}^{(n)}$ defined above provides a geometric model for the relative homology groups over $\mathbb{Q}$.  
This in turn yields a homological stability result, originally obtained by Nesterenko and Suslin~\cite{nes-sus1990}; for earlier stability results for general linear groups, see~\cite{vdkallen1980}.

\begin{theorem}\label{theorem4.7}
  For any \(n \ge 1\), there exists an isomorphism:
  \[
  H_\ast(\GL_n, \GL_{n-1}; \mathbb{Q}) \cong H_{\ast+1}\left((\widetilde{A}^{(n)}_{\ast})_{\GL_n}\right) \otimes \mathbb{Q}.
  \]
\end{theorem}

\begin{proof}
  Throughout the proof we will assume that the complex $\widetilde{A}_{\ast}^{(n)}$ is already over $\Q$.
  Consider the hyperhomology spectral sequence for the complex \(\widetilde{A}^{(n)}_*\):
  \[
  E_{p,q}^1 = H_p(\GL_n; \widetilde{A}^{(n)}_q) \implies H_{p+q}(\GL_n; \widetilde{A}^{(n)}_*) = 0,
  \]
  where the vanishing follows from the acyclicity of \(\widetilde{A}^{(n)}_*\).

  By the previous results, the \(E^1\)-page of this spectral sequence has the form:

\begin{center}
\begin{tikzcd}[cramped,sep=tiny]
	{} \\
	& {H_{p}(G_{n}; \mathbb{Q})} & {H_{p}(G_{n-1}; \mathbb{Q})} & 0 & \cdots & 0 & 0 & 0 \\
	& \vdots & \vdots & 0 &  & 0 & 0 & 0  \\
	& \vdots & \vdots & 0 &  & 0 & 0 & 0  \\
	& {H_{0}(G_{n}; \mathbb{Q})} & {H_{0}(G_{n-1}; \mathbb{Q})} & 0 & \cdots & 0  & {(\widetilde{A}_{n+1}^{(n)})_{G_{n}}} & {(\widetilde{A}_{n+2}^{(n)})_{G_{n}}} \\
	{} &&&&&&& {} & {}
	\arrow[color={rgb,255:red,92;green,92;blue,214}, from=2-3, to=2-2]
	\arrow[color={rgb,255:red,92;green,92;blue,214}, from=5-3, to=5-2]
	\arrow[color={rgb,255:red,92;green,92;blue,214}, from=5-8, to=5-7]
	\arrow[color={rgb,255:red,92;green,92;blue,214}, from=5-7, to=5-6]
	\arrow[shift right=3, from=6-1, to=1-1]
	\arrow[shift left, from=6-1, to=6-9]
\end{tikzcd}  
\end{center}

Define:
  \[
  \chi_i = \coker\left(H_i(\GL_{n-1}; \mathbb{Q}) \to H_i(\GL_n; \mathbb{Q})\right), \quad 
  \kappa_i = \ker\left(H_i(\GL_{n-1}; \mathbb{Q}) \to H_i(\GL_n; \mathbb{Q})\right).
  \]

  Since the spectral sequence converges to zero, we deduce:
  \[
  \chi_i = \kappa_i = 0 \quad \text{for } 0 \leq i \leq n - 2,
  \]
  and furthermore \(\chi_{n-1} = 0\), as all differentials into or out of these positions are zero.

  From the standard long exact sequence of relative homology:
  \begin{equation}\label{les_3}
      H_*(\GL_{n-1}; \mathbb{Q}) \to H_*(\GL_n; \mathbb{Q}) \to H_*(\GL_n, \GL_{n-1}; \mathbb{Q}),
  \end{equation}
  it follows that \(H_i(\GL_n, \GL_{n-1}; \mathbb{Q}) = 0\) for \(i \leq n - 2\), so we obtained an isomorphism for i-th homology with $i \leq n - 2$, because of \ref{lemma4.4}.  

  To see that \(H_{n-1}(\GL_n, \GL_{n-1}; \mathbb{Q}) = 0\), consider the segment:
  \[
  H_{n-1}(\GL_n; \mathbb{Q}) \xrightarrow{0} H_{n-1}(\GL_n, \GL_{n-1}; \mathbb{Q}) \xrightarrow{0} H_{n-2}(\GL_{n-1}; \mathbb{Q}) \xrightarrow{\sim} H_{n-2}(\GL_n; \mathbb{Q}).
  \]
  The first map vanishes since \(\chi_{n-1} = 0\), and the last map is an isomorphism.
  
  This forces \(H_{n-1}(\GL_n, \GL_{n-1}; \mathbb{Q}) = 0\), thus for $i \leq n - 1$, we have:
  \[
  H_{i}(\GL_n, \GL_{n-1}; \mathbb{Q}) \cong H_{i+1}\left((\widetilde{A}^{(n)}_{\ast})_{\GL_n}\right) = 0.
  \]

  The convergence of the spectral sequence yields short exact sequences for \(k \geq 0\):
  \[
  0 \to \chi_{n+k} \hookrightarrow H_{n+k+1}\left((\widetilde{A}^{(n)}_*)_{\GL_n}\right) \twoheadrightarrow \kappa_{n+k-1} \to 0.
  \]

  Consequently we obtain the following long exact sequence:
  \[
  \cdots \to H_{n+k}(\GL_n; \mathbb{Q}) \to H_{n+k+1}\left((\widetilde{A}^{(n)}_*)_{\GL_n}\right) \to H_{n+k-1}(\GL_{n-1}; \mathbb{Q}) \to H_{n+k-1}(\GL_n; \mathbb{Q}).
  \]

  Applying the five lemma to the long exact sequence (\ref{les_3}) implies:
  \[
   H_{n+k}(\GL_n, \GL_{n-1}; \mathbb{Q}) \cong H_{n+k+1}\left((\widetilde{A}^{(n)}_*)_{\GL_n}\right) \quad \text{for all } k \geq 0,
  \]
  that complets the proof.
\end{proof}

As a consequence of Theorem~\ref{theorem4.7} and the spectral sequence structure analyzed therein, we obtain:

\begin{cor}\label{cor:homology_stability}
    The standard embedding \(\GL_{n-1} \hookrightarrow \GL_n\) induces isomorphisms:
    \[
    H_{i}\left(\GL_{n-1}; \Z\left[\tfrac{1}{(n-1)!}\right]\right) \cong H_{i}\left(\GL_n; \Z\left[\tfrac{1}{(n-1)!}\right]\right) \quad \text{for all } i \leq n - 2,
    \]
    with the induced map being surjective for \(i = n - 1\).
\end{cor}

This result was originally established by Nesterenko and Suslin for  the \(\Z\) coefficients~\cite{nes-sus1990}.  
We provide an alternative and simplier proof here, for the localized coefficients \(\Z[\tfrac{1}{(n-1)!}]\). 

The proof of the remaining part~\ref{maincor2} is somewhat technical; however, we include it here because it provides a clear description of the long differentials that arise in the spectral sequence of \ref{theorem4.7}.

\section{Relation to Milnor's $K$-groups}

For this section, let us denote  
\[
t_i := \mathrm{ker}\left(\widetilde{A}_{i}^{(n)} \rightarrow \widetilde{A}_{i-1}^{(n)}\right) \cong \mathrm{coker}\left(\widetilde{A}_{i-1}^{(n)} \rightarrow \widetilde{A}_{i}^{(n)}\right).
\]

In particular, there is a short exact sequence  
\[
t_1 \hookrightarrow \widetilde{A}_{1}^{(n)} \twoheadrightarrow \widetilde{A}_{0}^{(n)}.
\]
From the corresponding long exact sequence, we obtain a map  
\[
\begin{tikzcd}[cramped]
	{H_{n}(\GL_{n})} && {H_{n-1}(\GL_{n}, t_1)} \\
	& {H_{n}(\GL_n) / H_{n}(\GL_{n-1})}
	\arrow["\partial", from=1-1, to=1-3]
	\arrow[from=1-1, to=2-2]
	\arrow["\partial"', hook, from=2-2, to=1-3]
\end{tikzcd}
\]

which factors through \(H_{n}(\GL_n, \widetilde{A}^{(n)}_1) \cong H_{n}(\GL_{n-1})\) (see~\ref{lemma4.2}) by exactness.  
Note that the top map \(\partial\) is the connecting homomorphism from the long exact sequence; as we will only use the factored map, let us also denote it by \(\partial\).

Furthermore, by standard dimension shifting we obtain a map which, after passing to coefficients with \((n-1)!\) inverted, becomes an isomorphism:  

\[
H_{n-1}(\GL_n, t_1) \longrightarrow H_{n-2}(\GL_n, t_2) \longrightarrow \cdots \longrightarrow H_{0}(\GL_n, t_n).
\]

Indeed, each of the maps above is a connecting homomorphism arising from the corresponding long exact sequence associated to the short exact sequence  
\[
t_{i} \hookrightarrow \widetilde{A}_{i}^{(n)} \twoheadrightarrow t_{i-1}.
\]
Since the homology groups \(H_{\ast}(\GL_n, \widetilde{A}_{i}^{(n)})\) vanish in the relevant range after inverting \((n-1)!\) in the coefficients, we obtain the desired isomorphism.

We have thus established the following:

\begin{prop}
  There exists a sequence of maps induced by connecting homomorphisms
  \[
  H_{n-1}(\GL_n; t_1) \xrightarrow{\sim} H_{n-2}(\GL_n; t_2) \xrightarrow{\sim} \cdots \xrightarrow{\sim} H_{0}(\GL_n; t_n),
  \] 
    which form isomorphisms after inverting \((n-1)!\) in the coefficients.
\end{prop}

\begin{convention}
From this point on, we assume that \((n-1)!\) is inverted in the entire complex \(\widetilde{A}_{\ast}^{(n)}\).
\end{convention}

Recall that we considered a spectral sequence built on the complex \(\widetilde{A}_{\ast}^{(n)}\); we now focus on the second page \(E_{n+1, 0}^{2}\).  
By previous calculations~\ref{lemma4.4} we have the following description \(E_{n, 0}^{1} \cong 0\) and \(E_{n+i, 0}^{1} \cong (\widetilde{A}_{n+i}^{(n)})_{G_n}\) for \(i \geq 1\). Passing to the second page, we obtain  

\[
E_{n+1, 0}^{2} \cong \mathrm{coker}\left(H_{0}(\GL_n, \widetilde{A}_{n+2}^{(n)}) \longrightarrow H_{0}(\GL_n, \widetilde{A}_{n+1}^{(n)})\right).
\]

We now prove the following proposition:  

\begin{prop}
  \(H_{0}(\GL_n, t_n) \cong E_{n+1,0}^{2}\).
\end{prop}

\begin{proof}
  From the discussion above, we have
  \[
  E_{n+1, 0}^{2} \cong \mathrm{coker}\left(H_{0}(\GL_n, \widetilde{A}_{n+2}^{(n)}) \longrightarrow H_{0}(\GL_n, \widetilde{A}_{n+1}^{(n)})\right),
  \]
  while the short exact sequence
  \[
  t_{n+1} \hookrightarrow \widetilde{A}_{n+1}^{(n)} \twoheadrightarrow t_{n}
  \]
  yields
  \[
  H_{0}(\GL_n, t_n) \cong \mathrm{coker}\left(H_{0}(\GL_n, t_{n+1}) \longrightarrow H_{0}(\GL_n, \widetilde{A}_{n+1}^{(n)})\right).
  \]
  Since the tensor product preserves surjections, the induced map
  \[
  H_{0}(\GL_n, \widetilde{A}_{n+2}^{(n)}) \longrightarrow H_{0}(\GL_n, t_{n+1})
  \]
  is surjective. Therefore, the two cokernels above coincide, establishing the claim.
\end{proof}

The general position conditions ensure that the zeroth homology group \(H_{0}(\GL_n, \widetilde{A}_{n+1}^{(n)})\) is generated by \(n \times (n+1)\) matrices of the form \((e_1, \ldots, e_n, \bar{\alpha})\),  
where \(e_i\) denote the standard basis vectors in \(F^{n}\), and \(\bar{\alpha}\) is an arbitrary element of \(F^{n}\) with all nonzero coordinates.

We denote such a generator by \([\bar{\alpha}] = \{ \alpha_1, \ldots, \alpha_n \}\), recording only the data on which the generator depends.  
Since \(\widetilde{A}_{n+1}^{(n)}\) is an abelian group generated by \((n+1)\)-tuples of vectors satisfying the general position condition, its coinvariants form an abelian group generated by elements of the form \(\{\alpha_1, \ldots, \alpha_n\}\), where each \(\alpha_i \in F^{\ast}\) is arbitrary.

To determine the structure of \(E_{n+1, 0}^{2}\) in terms of generators and relations, we describe which relations arise from the cokernel of the map  
\[
H_{0}(\GL_{n}, \widetilde{A}_{n+2}^{(n)}) \longrightarrow H_{0}(\GL_{n}, \widetilde{A}_{n+1}^{(n)}).
\]

\begin{prop}
   The group \(E_{n+1, 0}^{2}\) is generated by the elements
   \[
   \{\alpha_1, \ldots, \alpha_n\} := (e_1, \ldots, e_n, \bar{\alpha}),
   \]
   subject to the relations
   \[
   \sum_{i=1}^{n} (-1)^{i+n} 
   \{\alpha_1(\beta_1-\beta_i), \ldots, \widehat{\alpha_i(\beta_i - \beta_i)}, \ldots, \alpha_n(\beta_n - \beta_i), \beta_i\}
   = \{\alpha_1 \beta_1, \ldots, \alpha_n \beta_n\} - \{\alpha_1, \ldots, \alpha_n\},
   \] 
   and 
   \[
  \{\alpha_1, \ldots, \alpha_n\} = (-1)^{\sign(\sigma)}\{\alpha_{\sigma^{-1}(1)}, \ldots, \alpha_{\sigma^{-1}(n)}\},
   \]
   where \(\alpha_i\) and \(\beta_i\) are such that the vectors  
   \[
   \left( e_1, \ldots, e_n, (\alpha_1, \ldots, \alpha_n)^{t}, (\alpha_1 \beta_1, \ldots, \alpha_n \beta_n)^{t} \right)
   \]
   together satisfy the general position condition.
\end{prop}

\begin{proof}
  Consider an arbitrary generator of \((\widetilde{A}_{n+2}^{(n)})_{\GL_{n}}\).  
  Up to multiplication by an element of \(\GL_n\), we may assume it has the form \((e_1, \ldots, e_n, \bar{\alpha}, \bar{\beta})\),  
  where, as usual, the \(e_i\) are the vectors of the standard basis and \(\bar{\alpha}, \bar{\beta}\) satisfy the general position condition.  
  For convenience, set
  \([\bar{\alpha}] = \{\alpha_1, \ldots, \alpha_n\}\) and \([\bar{\beta}] = \{\alpha_1 \beta_1, \ldots, \alpha_n \beta_n\}\).  
  It is clear that the second vector can be expressed in this form.

  All relations in \(E_{n+1, 0}^{2}\) arise as linear combinations of the images of such generators under the standard differential.  
  Hence, we have the relation
  \[
  \sum_{i=1}^{n} (-1)^{i-1} (e_1, \ldots, \widehat{e_i}, \ldots, e_n, \bar{\alpha}, \bar{\beta})
  = (-1)^{n+1}(e_1, \ldots, e_n, \bar{\beta})
    + (-1)^{n+2}(e_1, \ldots, e_n, \bar{\alpha}).
  \]

  A straightforward computation shows that
  \[
  (e_1, \ldots, \widehat{e_i}, \ldots, e_n, \bar{\alpha}, \bar{\beta})
  = \{\alpha_1(\beta_1-\beta_i), \ldots,
  \widehat{\alpha_i(\beta_i - \beta_i)}, \ldots,
  \alpha_n(\beta_n - \beta_i), \beta_i\},
  \]
  and thus, in terms of the generators, we obtain
  \[
   \sum_{i=1}^{n} (-1)^{i+n} 
   \{\alpha_1(\beta_1-\beta_i), \ldots, \widehat{\alpha_i(\beta_i - \beta_i)}, \ldots, \alpha_n(\beta_n - \beta_i), \beta_i\}
   = \{\alpha_1 \beta_1, \ldots, \alpha_n \beta_n\} - \{\alpha_1, \ldots, \alpha_n\}.
  \]
\end{proof}

\begin{remark}
The conditions on \(\alpha_i\) and \(\beta_i\) can be reformulated as follows.  
The general position of any \(n-1\) vectors from the standard basis together with one vector from either \(\bar{\alpha}\) or \(\bar{\beta}\) is equivalent to requiring that all \(\alpha_i\) and \(\beta_i\) are nonzero.  
The general position of any \(n-2\) vectors from the standard basis together with both \(\bar{\alpha}\) and \(\bar{\beta}\) is equivalent to requiring that all \(\beta_i\) are distinct.
\end{remark}

Remarkably, the relations in \(E^{2}_{n+1, 0}\) also hold in the \(n\)-th Milnor \(K\)-theory; see Proposition~3.4 in~\cite{nes-sus1990}.  
Hence, there is an induced map
\[
E_{n+1, 0}^{2} \longrightarrow K^{M}_{n}(F),
\]
which sends each generator \(\{\alpha_1, \ldots, \alpha_n\}\) to the corresponding generator of Milnor \(K\)-theory, denoted by the same symbol.  
Consequently, we obtain the following sequence of maps:
\[
H_{n}(\GL_{n}(F))/H_{n}(\GL_{n-1}(F))
\xrightarrow{\partial}
H_{n-1}(\GL_{n}, t_1)
\xrightarrow{\sim}
E_{n+1, 0}^{2}
\longrightarrow
K^{M}_{n}(F).
\]

There is a natural map from Milnor \(K\)-theory \(K^{M}_{n}(F)\) to \(H_{n}(\GL_n(F))/H_{n}(\GL_{n-1}(F))\), constructed via the cross product in homology.  
Indeed, there is a cross-product map
\[
\underbrace{F^{\ast} \otimes \cdots \otimes F^{\ast}}_{n}
\cong
H_{1}(\GL_{1}) \otimes \cdots \otimes H_{1}(\GL_{1})
\xrightarrow{\times}
H_{n}(\GL_{n}(F)).
\]  
The fact that the cross-product map, composed with the quotient by \(H_{n}(\GL_{n-1})\), factors through Milnor \(K\)-theory was proved in another paper by Suslin; see~\cite{sus1985}, Corollary~2.7.2.

Hence, we obtain the following diagram: 

\[\begin{tikzcd}\label{main_diagram}
	& {\ker(H_{n}(\GL_{n-1}) \rightarrow H_{n}(\GL_{n}))} \\
	{H_{n}(\GL_{n})/H_{n}(\GL_{n-1})} & {H_{n-1}(\GL_{n}, t_1)} && {E^{2}_{n+1, 0}} & {K^{M}_{n}(F)}
	\arrow["\partial", hook, from=2-1, to=2-2]
	\arrow[two heads, from=2-2, to=1-2]
	\arrow["\sim", from=2-2, to=2-4]
	\arrow[two heads, from=2-4, to=2-5]
	\arrow[shift left, curve={height=-30pt}, from=2-5, to=2-1]
\end{tikzcd}\]

Our goal is to prove that all horizontal morphisms in the diagram are isomorphisms, while the vertical map on the left is the zero map.  
We begin by showing that the composition starting in Milnor \(K\)-theory and ending in it is the identity map.

Let us now describe the morphisms in the diagram above more explicitly, starting from the generators of \(E_{n+1, 0}^{2}\).  
Each generator \(\{\alpha_1, \ldots, \alpha_n\}\) maps to the corresponding generator \(\{\alpha_1, \ldots, \alpha_n\}\) in the Milnor \(K\)-group \(K^{M}_{n}(F)\).

For any \(n\)-tuple \((g_1, g_2, \ldots, g_n)\) of pairwise commuting elements in \(\GL_{n}(F)\), consider the homology class
\[
c(g_1, g_2, \ldots, g_n) \in H_n(\GL_{n})
\]
represented by the cycle
\[
\sum_{\sigma \in S_n} \operatorname{sign}(\sigma)
\,[g_{\sigma(1)} \mid g_{\sigma(2)} \mid \cdots \mid g_{\sigma(n)}] \otimes 1
\]

\begin{prop}
The cross product map from \(K^{M}_{n}(F)\) to \(H_{n}(\GL_{n}(F)) / H_{n}(\GL_{n-1}(F))\) has the explicit description (see~\cite{guin1989}, \cite{sus1985}):
\[
\{\alpha_1, \ldots, \alpha_n\}
\longmapsto
c(D_{1}(\alpha_1), \ldots, D_{n}(\alpha_n))
\,[H_{n}(\GL_{n-1})],
\]
where \(D_{i}(a)\) is the diagonal matrix with \(a\) in the \(i\)-th position and \(1\) elsewhere.
\end{prop}

We now examine the action of the remaining horizontal maps in diagram \ref{main_diagram} on the generator \(c(D_{1}(\alpha_1), \ldots, D_{n}(\alpha_n))\) via the standard free resolution. We begin with the boundary map \(\partial: H_{n}(\GL_{n})/H_{n}(\GL_{n-1}) \to H_{n-1}(\GL_{n}, t_1)\).

For clarity, we omit the \(\alpha\)-indices as they are uniquely determined by the indices of \(D\). The connecting morphism \(\partial\) is constructed from the long exact sequence, as illustrated by the diagram below:
\[\begin{tikzcd}[column sep=scriptsize,row sep=scriptsize]
	& {c(D_1, \ldots, D_n)\otimes e_n} & {c(D_1, \ldots, D_n)\otimes1} \\
	{\partial(c(D_1, \ldots, D_n) \otimes 1)} & {\partial(c(D_1, \ldots, D_n))\otimes e_n \in{H_{n-1}(\GL_n, \widetilde{A}_1)}}
	\arrow[from=1-2, to=1-3]
	\arrow[from=1-2, to=2-2]
	\arrow[dashed, from=1-3, to=2-1]
	\arrow[from=2-1, to=2-2]
\end{tikzcd}\]

The construction lifts the class \(c(D_{1}(\alpha_1), \ldots, D_{n}(\alpha_n)) \otimes 1\) to \(c(D_{1}(\alpha_1), \ldots, D_{n}(\alpha_n)) \otimes e_n\) in \(H_{n}(\GL_n, \widetilde{A}_1)\), then applies the differential in the free resolution. 

We can explicitly describe the image of this class under the boundary map.

\begin{prop}\label{description of the first partial}
With notation as above, the boundary of \(c(D_{1}(\alpha_1), \ldots, D_{n}(\alpha_n))\) satisfies
\[
\partial\!\bigl( \, c(D_{1}(\alpha_1), \ldots, D_{n}(\alpha_n)) \;\otimes\; 1 \, \bigr)
\;=\;
c(D_{1}(\alpha_1), \ldots, D_{n-1}(\alpha_{n-1}))
\;\otimes\; 
(-1)^{n} \, \bigl( \, e_n - D_{n}(\alpha_{n})e_n \, \bigr) \in H_{n-1}(\GL_n, t_1).
\]
\end{prop}
\begin{proof}
First, consider the inner face operators \(\delta_{i}\) for \(0 < i < n\) applied to the cycle. 
\[
\begin{aligned}
\delta_{i}\!\left(c(D_{1}(\alpha_1), \ldots, D_{n}(\alpha_n))\right) \otimes e_n
	&= \delta_{i} \!\left(
	\sum_{\sigma \in S_n} \operatorname{sign}(\sigma)
\,[D_{\sigma(1)}(\alpha_{\sigma(1)}) \mid D_{\sigma(2)}(\alpha_{\sigma(2)}) \mid \cdots \mid D_{\sigma(n)}(\alpha_{\sigma(n)})]
	\right) \otimes e_n \\
\end{aligned}
\]

We group permutations that differ only by transposing \(i\) and \(i+1\). The corresponding summands then have opposite signs and cancel:
\[
\begin{aligned}
\delta_{i}\Big(
[D_{\sigma(1)} \mid \ldots \mid 
   D_{\sigma(i)} \mid D_{\sigma(i+1)} \mid \ldots \mid D_{\sigma(n)}] - [D_{\sigma(1)} \mid \ldots \mid D_{\sigma(i+1)} \mid D_{\sigma(i)} \mid  \ldots \mid D_{\sigma(n)}]
\Big) = 0.
\end{aligned}
\]

Therefore, the only terms that contribute are those corresponding to \(\delta_{0}\) and \(\delta_{n}\).

We now group summands differing by the transposition \((1, n)\) and examine their sum:
\[
\begin{aligned}
&\delta_{0}\!\Big(
\operatorname{sign}(\sigma)\,
[D_{\sigma(1)} \mid D_{\sigma(2)} \mid \ldots \mid D_{\sigma(n)}]
\Big)\otimes e_n \\[2pt]
&\quad +\;
\delta_{n}\!\Big(
(-1)^{n-1}\operatorname{sign}(\sigma)\,
[D_{\sigma(2)} \mid D_{\sigma(3)} \mid \ldots \mid D_{\sigma(n)} \mid D_{\sigma(1)}]
\Big)\otimes e_n.
\end{aligned}
\]

This sum vanishes unless \(\sigma(1) = n\), since the action on \(e_n\) is trivial otherwise. 
Summing over all \(\sigma\) with \(\sigma(1) = n\) and applying the action on \(e_n\) completes the proof.
\end{proof}

By the same reasoning, we obtain the following general statement.  

\begin{prop}\label{description of partial}
Suppose \(x \in t_i\) with preimage \(\tilde{x} \in \widetilde{A}_{i+1}\) such that \(D_{j}(\alpha_{j})\) acts trivially on \(\tilde{x}\) for \(j < n - i\).
Then the following formula holds:  
\[
\partial\!\bigl(c(D_{1}(\alpha_1), \ldots, D_{n-i}(\alpha_{n-i}))\bigr) \otimes x
\;=\;
c(D_{1}(\alpha_1), \ldots, D_{n-i-1}(\alpha_{n-i-1}))
\otimes 
(-1)^{n-i}\!\bigl(\tilde{x} - D_{n-i}(\alpha_{n-i})\tilde{x}\bigr).
\]
\end{prop}

We now explicitly describe the composite map:
\[
K^{M}_{n}(F)
\longrightarrow
H_{n}(\GL_n) / H_{n}(\GL_{n-1})
\xhookrightarrow{\;\partial\;}
H_{n-1}(\GL_{n}, t_1)
\xrightarrow{\;\cong\;}
E_{n+1,0}^{2}.
\]
\todo[inline]{Rewrite the proof}
\begin{prop}
The composition sends each generator \(\{\alpha_1, \ldots, \alpha_n\}\) of \(K^M_n(F)\) to
\[
\sum_{\{i_1, \ldots, i_j\} \subset \{1, \ldots, n\}}
(-1)^{j}\,
\{\alpha_1, \ldots, 1, \ldots, 1, \ldots, \alpha_n\},
\]
where each summand has \(1\) at positions \(i_1, \ldots, i_j\) and the corresponding \(\alpha_i\) in the remaining positions.
\end{prop}

\begin{proof}
We construct sequences \((x_i)_{i=1}^{n}\) and \((\tilde{x}_i)_{i=0}^{n}\) satisfying the following conditions:
\begin{itemize}
  \item \(x_i \in t_i\),
  \item \(\tilde{x}_i \in \widetilde{A}_{i + 1}\),
  \item \(\tilde{x}_i \mapsto x_i \text{ under the map } A_{i + 1} \longrightarrow t_i.\)
\end{itemize}
We choose these sequences so that \(D_{j}(\alpha_{j})\) acts trivially on each \(\tilde{x}_i\) for \(j < n - i\). These choices are reflected in the following commutative diagrams for each \(i = 0, \ldots, n-1\):                  
\[
\begin{tikzcd}[cramped,column sep=scriptsize,row sep=scriptsize]
	&& {c(D_{1}, \ldots, D_{n-i}) \otimes \tilde{x}_i} && {c(D_{1}, \ldots, D_{n-i}) \otimes x_i} \\
	\\
	{\partial(c(D_{1}, \ldots, D_{n-i}) \otimes x_{i})} && {\partial\left(c(D_{1}, \ldots, D_{n-i}) \right) \otimes \tilde{x}_{i}}
	\arrow[from=1-3, to=1-5]
	\arrow[from=1-3, to=3-3]
	\arrow[dashed, from=1-5, to=3-1]
	\arrow[from=3-1, to=3-3]
\end{tikzcd}
\]
We define these sequences inductively by setting
\[
x_0 = 1 \qquad \text{and} \qquad \tilde{x}_0 = e_n,
\]
as in Proposition~\ref{description of the first partial}, then by Proposition~\ref{description of partial}:
\[
x_i = (-1)^{\,n-i+1} \,\bigl(\tilde{x}_{i-1} - D_{\,n-i+1}(\alpha_{n-i+1})\, \tilde{x}_{i-1}\bigr).
\]
For the inductive step, we define each \(\tilde{x}_{i}\) by choosing a preimage of \(x_{i}\) that preserves the required properties.
Specifically, merging each summand with \(e_{n-i}\) on the left yields such a preimage (by the same reasoning as in Proposition~\ref{prop2.3}), ensuring that \(D_{j}(\alpha_{j})\) acts trivially on \(\tilde{x}_i\) for \(j < n - i\). 

Finally, to track the image in \(E_{n+1, 0}^{2}\), we lift \(x_{n}\) via the surjection \((A_{n+1})_{\GL_{n}} \twoheadrightarrow (t_n)_{\GL_{n}}\) by adjoining the column vector \((\alpha_1, \ldots, \alpha_n)^{t}\) to each summand. This construction ensures all columns satisfy the general position condition, and expressing each summand as a generator establishes the claim.
\end{proof}

\begin{cor}
The composition map
\[
K^{M}_{n}(F) \rightarrow H_{n}(\GL_{n})/H_{n}(\GL_{n-1}) \hookrightarrow H_{n-1}(\GL_{n}, t_1) \xrightarrow{\sim} E_{n+1, 0}^{2} \rightarrow K^{M}_{n}(F)
\]
is the identity.
\end{cor}

\begin{proof}
The corollary follows from the previous proposition combined with the fact that in Milnor \(K\)-theory, the generator \(\{\alpha_1, \ldots, \alpha_n\}\) vanishes whenever any entry equals \(1\).
\end{proof}

Moreover, the same relation—that \(\{\alpha_1, \ldots, \alpha_n\} = 1\) whenever some \(\alpha_i = 1\)—also holds in \(E_{n+1, 0}^{2}\).  
Although conceptually straightforward, the proof requires detailed calculations in \(E_{n+1, 0}^{2}\), which we omit for brevity.  
As a consequence, the composition 
\[
E_{n+1, 0}^{2} \rightarrow K^{M}_{n}(F) \rightarrow H_{n}(\GL_{n})/H_{n}(\GL_{n-1}) \hookrightarrow H_{n-1}(\GL_{n}, t_1) \xrightarrow{\sim} E_{n+1, 0}^{2}
\]
is also the identity. 

These results allow us to recover the two main theorems of Nesterenko and Suslin~\cite{nes-sus1990}.

\begin{cor}
The map induced by the inclusion \(GL_{n-1} \hookrightarrow \GL_{n}\) is injective on \((n-1)\)-st homology after inverting \((n-1)!\).  
Moreover, the cokernel of the induced map on the \(n\)-th homology is isomorphic to the \(n\)-th Milnor \(K\)-group:
\[
H_{n-1}(\GL_{n-1}, \Z[\tfrac{1}{(n-1)!}]) \hookrightarrow H_{n-1}(\GL_{n}, \Z[\tfrac{1}{(n-1)!}]),
\]
and
\[
H_{n}(\GL_{n})/H_{n}(\GL_{n-1}) \otimes \Z[\tfrac{1}{(n-1)!}] \cong K^{M}_{n}(F).
\]
\end{cor}

\begin{remark}
Nesterenko and Suslin originally proved these results integrally, employing a more elaborate machinery.  
We believe, however, that the approach developed in this paper not only simplifies their arguments but also connects naturally to the recent results of Galatius–Kupers–Randal-Williams~\cite{G-K-R2020}, who established that the map
\(H_{n}(\GL_{n-1}, k) \to H_{n}(\GL_{n}, k)\)
is injective whenever \((n-1)!\) is invertible in the coefficient field \(k\).
\end{remark}

\todo[inline]{Add a general statement with a description of long differentials}

\todo[inline]{Try to prove the remaining cases of Suslin's Injectivity conjecture via this description}

\todo[inline]{Section about multiplication?}

\todo[inline]{Section about projective complexes?}

\bibliographystyle{plain}
\bibliography{ref}

\end{document}